\documentstyle[12pt,psfig]{article}
\begin{document}
\title{TQFT and Whitehead's manifold}
\author{Louis FUNAR\\
Institut Fourier, BP 74, Univ.Grenoble I, Math\'ematiques \\
38402 Saint-Martin-d'H\`eres cedex}
\maketitle
\abstract{The aim of this note is to derive some invariants at
infinity for open 3-manifolds in the framework of Topological
Quantum Field Theories.
These invariants may be used to test if an open manifold is
simply connected at infinity as we done for
Whitehead's manifold in case of the $sl_{2}({\bf C})$-TQFT in level 4.}

\section{Introduction}
The aim of this paper is to introduce the invariants at infinity for open
 3-manifolds deduced from
Topological Quantum Field Theories (abbrev. TQFT). We are able to compute the
invariant of the classical Whitehead manifold associated to the simplest
non-abelian TQFT based on $sl_2({\bf C})$ in level 4 and find it is not
trivial.
As a consequence the manifold is not simply connected at infinity. Further
 developpements
accreditating the idea that TQFT may share some light on the topology of
open 3-manifold
are pursued in a further paper.

In the first section we discuss the general TQFT invariants at infinity
 $Z_{\infty}(W)$ for an open
3-manifold $W$. In order to preserve the self-contained character of this paper
 we outline the
definition of link and 3-manifold invariants of Witten and Reshetikhin-Turaev
\cite{Wi,Re-Tu} based on the quantum $sl_2({\bf C})$ at roots of unity.
According to a general result each multiplicative invariant for closed
3-manifolds
extends canonically to a TQFT (\cite{Fun5}). We follow the surgical approach
 used by
Turaev in \cite{Tu2} to get the explicit description of the TQFT for
cobordisms.

In the last section we use these results to compute effectively the TQFT
invariant at infinity for the Whitehead manifold. Recall (\cite{Wh}) this is
defined
 as
follows: Let $T_1\hookrightarrow T_0$ be the embedding of the solid tori from
 picture 1.
\begin{figure}
\centering
\hspace{1cm}
\psfig{figure=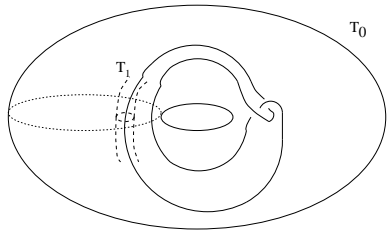,width=5cm}
\caption{The inclusion of tori}
\end{figure}
There exists an homeomorphism $h$ of $S^3$ so that $h(T_1)=T_0$. Consider
 the open manifold
$Wh=\bigcup_{n\geq 0}h^n(T_0)$. Then $Wh$ is the typical example of a
contractible
open 3-manifold which is not homeomorphic to ${\bf R}^3$. The precise reason
 is that
$Wh$ is not simply connected at infinity (i.e. not every compact may be
engulfed
 in
a compact simply-connected submanifold).

Our main result states as follows:
\newtheorem{ttc}{Theorem}[section]
\begin{ttc}
The space $Z_{\infty}(W)\cong {\bf C}^3$ if $Z$ is the $sl_2({\bf C})$-TQFT
at level 4.
\end{ttc}
Since any open 3-manifold $W$ which is simply-connected at infinity must
 satisfy
$\dim Z_{\infty}(W)\leq 1$ for all reduced TQFT (see Proposition 2.2) we get
another proof for the non simple connectedness at infinity of $Wh$ .
It seems that $Z_{\infty}$ is a good test for the simple connectedness at
infinity
even if not enough
strong to distinguish among various non-homeomorphic  open 3-manifolds.

{\bf Acknowledgements}. This paper is based on author's PhD thesis at
University of
Paris-Sud, 1994. I'am grateful to my advisor Valentin Po\'enaru for proposing
 me this
problem and for many discussions we had on this subject, to Pierre Vogel and
Vladimir Turaev for their helpful comments and suggestions and to the referee
 for his
careful reading of the manuscript leading to a considerably improved version.

\section{Invariants at infinity from TQFT}
(2.1) {\bf On TQFT}. Recall \cite{At} that a TQFT in dimension 3 is a functor
$Z$
from
the category of oriented cobordisms into that of hermitian vector spaces. This
 means that
to a compact surface $S$ we associate an hermitian vector space $Z(S)$
depending only on
the topological type of $S$. The {\it quantum} character of the theory is
reflected
 in the
rules

$Z(\cup_i S_i)=\bigotimes_i Z(S_i)$, $Z(\emptyset)={\bf C}$, \\
which make the difference with the usual functors encountered in the algebraic
 topology.

Furthermore to an oriented cobordism $M$ so that $\partial M$ is split into two
disjoint
manifold $S$ and $T$ (the {\it incoming} and the {\it outgoing} boundaries
respectively, which are not necessary those given by the orientation)
we have assigned a morphism $Z(M):Z(S)\longrightarrow Z(T)$ satisfying the
natural compatibility relations between composition of morphisms and
cobordisms.
This is usually called an anomaly-free TQFT (see \cite{Tu2}). The main examples
yet constructed
have an anomaly from a certain group of roots of unity $\Gamma\subset U(1)$.
This means that the invariant associated to the composition of cobordisms
 $M$ and $N$ may be expressed as

$Z(M\circ N)=\gamma Z(M)\circ Z(N)$, with $\gamma\in \Gamma$. \\
 The usual way to deal with this ambiguity is to work with framed 3-manifolds
\cite{At1} or $p_1$-structures \cite{BHMV}. However the presence of an anomaly
 will be
irrelevant for the construction of invariants at infinity.

The examples we consider in this paper are {\it reduced} TQFT, namely they
satisfy the additional condition

$Z(S^2)\cong {\bf C}$.

All the TQFT from
quantum groups or quasi-quantum groups are reduced and we may restrict
 ourselves to
the study of reduced TQFT by the results of \cite{Fun5}.

(2.2) {\bf Open 3-manifolds}. Consider first $Z$ is an anomaly-free TQFT. Let
$W$ be an open 3-manifold without boundary. We choose an ascending
 sequence
of submanifolds $\{K_n\}$ fulfilling

$K_n\subset int(K_{n+1})$, $W=\cup_nK_n$. \\
Then $V_i=cl(K_{i+1}-K_i)$ are oriented cobordism from $\partial K_i$ to
$\partial K_{i+1}$. Here $int$ and $cl$ state for the interior and the closure
respectively. We get a sequence of linear maps

$Z(V_i):Z(\partial K_i)\longrightarrow Z(\partial K_{i+1})$, \\
which represent an inductive system of vector spaces. We define $Z_{\infty}(W)$
 be simply the
inductive limit of this system.
\newtheorem{brc}{Definition-Lemma}[section]
\begin{brc}
The vector space $Z_{\infty}(W)$ is the topological invariant at infinity
associated to the TQFT functor $Z$ and the open 3-manifold $W$.
\end{brc}
In fact it is simply to check the independence of $Z_{\infty}(W)$ on the choice
of
 the
exhaustion or the parametrizations of the intermediary boundaries $\partial
K_i$.
$\Box$

Remark that $Z_{\infty}(W)$ depends only on the structure at infinity of $W$:
if
$W'$ is another manifold so that $W$ and $W'$ are homeomorphic outside some
 compacts
then the associated spaces $Z_{\infty}(W)$ and $Z_{\infty}(W')$ are isomorphic.

Also  if $Z$ is a TQFT with anomaly this time, then the maps $Z(V_i)$ are
defined
up to the multiplication by some scalar from $\Gamma$. Nevertheless the space
$Z_{\infty}(W)$ is well determined.
It is only the hermitian structure which is lost when we pass to the limit
 $Z_{\infty}(W)$.

(2.3) {\bf h-1-connected manifolds}. The open 3-manifold $W$ is
{\it h-1-connected at infinity}
if each compact $K\subset W$ may be engulfed in a compact submanifold
 $Y\subset W$
with $H_1(Y)=0$. This is a condition slightly weaker than the simple
connectedness at
infinity.
\newtheorem{brc1}[brc]{Proposition}
\begin{brc1}
If $W$ is h-1-connected at infinity then $\dim Z_{\infty}(W)\leq 1$ for any
reduced TQFT.
\end{brc1}
Proof: It suffices to observe that a compact 3-manifold $Y$ with $H_1(Y)=0$ has
 the boundary
$\partial Y$ an union of spheres $S^2$, from an Euler characteristic argument.
If $\{K_n\}$ is an exhaustion of $W$ like in the introduction then there exists
compact submanifold $Y_n$ with $H_1(Y_n)=0$ and a function $r(n)>n$ so that

$K_n\subset int(Y_n)$ and $Y_n\subset int(K_{r(n)})$. \\
Then the map $Z(cl(K_{r(n)}-K_n))$ factors through $Z(\partial Y_n)\cong {\bf
C}$
 (because the TQFT is reduced)
hence the rank of the limit is at most 1. $\Box$

(2.4) {\bf The Whitehead manifold}. Consider now $Z$ be the level 4
 $sl_2({\bf C})$-TQFT of
Witten and Reshetikhin-Turaev (see the next section for complete definitions).
The main result of this paper stated in introduction asserts that
 $Z_{\infty}(Wh)\cong {\bf C}^3$.

We defer for the proof in section 4. We already notice that $Wh$ has periodic
ends and we may choose
$K_n=\bigcup_{0\leq j\leq n}h^j(T_0)$. Then all intermediary cobordisms $V_n$
are homeomorphic to $X=cl(T_0-T_1)$. It suffices therefore to compute
the linear map $Z(X):Z(\partial T_1)\longrightarrow Z(\partial T_0)$. For the
TQFT
 we are
working with, the space associated to a torus $Z(S^1\times S^1)$ is ${\bf C}^3$
(see the further section).
So the statement of the theorem is equivalent to the non-degeneracy of the
linear map $Z(X)$.

\section{The $sl_2({\bf C})$-TQFT}
We fix some integer $r>1$ called the level of the theory. The description we
outline follows from
\cite{Re-Tu,Ki-Me,Tu2}.

(3.1) {\bf Framed tangles}. Recall that a tangle $T$ is a 1-manifold properly
embedded
in the unit cube $I^3$ in ${\bf R}^3$ with
$\partial T\subset \{\frac{1}{2}\}\times I\times \partial I$, considered
up to isotopy rel boundary. If $\partial_-T=T\cap I^2\times \{0\}$ and
$\partial_+T=T\cap I^2\times \{1\}$ then $T$ is a {\it $(m,n)$-tangle} provided
that
$\mid \partial_-T\mid =m$ and $\mid \partial_+T\mid =n$. Thus a link is a
 (0,0)-tangle and a general
tangle consists of a link with a collection of proper arcs. We assume the
tangles are oriented and transverse to $I^2\times \partial I$.

A {\it framed tangle} is a tangle equipped with a framing of its normal bundle
(up
 to isotopy)
which is standard on the boundary. It is equivalent to a ribbon tangle from
\cite{Re-Tu} if we think the tangle is thickened to a ribbon in the direction
of the
second vector of the framing. Anyway framings may be specified by integers
assigned to the components of $T$.

Also one studies tangles using generic projections, called {\it diagrams}, onto
$\{0\}\times I^2$ having only ordinary double points. We assume
the framing considered
is the blackboard one, in which the second vector is parallel to
$\{0\}\times I^2$, and further
coincides with the 0-framing by eventually adding the
necessary number of kinks.

It is simply to check that every tangle diagram may
be factored into the elementary tangles
from picture 2.
\begin{figure}
\centering
\hspace{1cm}
\psfig{figure=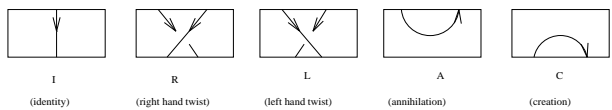,width=8cm}
\caption{Elementary tangles}
\end{figure}
There are well-known Reidemester moves describing the
local moves necessary and sufficient to obtain two
framed tangle diagrams one from the other if they are
coming from the same framed tangle.
For the sake of completeness we pictured them in
figure 3 (see \cite{Re-Tu}). Notice the orientations are
arbitrary and the framing is the blackboard one.
\begin{figure}
\centering
\hspace{1cm}
\psfig{figure=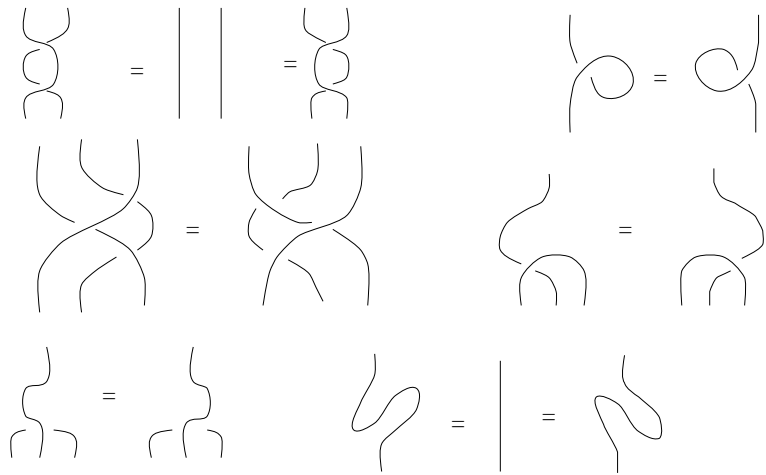,width=10cm}
\caption{Reidemester moves for tangles}
\end{figure}

(3.2) {\bf Quasi-Triangular Hopf algebras}. We discuss the quantum
$sl_2({\bf C})$ which is the main
example of a ribbon Hopf algebra. Recall $sl_2({\bf C})$ is
3-dimensional as vector space
and the Lie bracket is given (in terms of preferred generators) by:

$[H,X]=2X$, $[H,Y]=-2Y$, $[X,Y]=H$. \\
The {\it universal envelopping algebra} $U=U(sl_2({\bf C}))$ is the associative
algebra
over ${\bf C}$ generated by $X,Y,H$ and the relations from above. Notice that
there exists an unique
$k$-dimensional irreducible $sl_2({\bf C})$-module $V^k$, for each integer $k$,
which has also an $U$-module structure. Also $U$ is a Hopf algebra when
 endowed
with the comultiplication $\Delta:U\longrightarrow U\otimes U$ given by
$\Delta(u)=u\otimes 1 +1\otimes u$, antipode $s:U\longrightarrow U$ given by
$s(u)=-u$, and counit $\varepsilon:U\longrightarrow {\bf C}$ determined by
$\varepsilon(u)=0$, for all $u$ which are Lie polynomials in $X,Y,H$.

Now the {\it quantized universal envelopping algebra} $U_h=U_h(sl_2{\bf C}))$
is
defined
as $U[[h]]$ (the formal series in $h$) with the same relations as $U$ excepting
for
$[X,Y]=H$ which is replaced by
\begin{equation}
[X,Y]=[H]=\frac{e^{\frac{h}{2}H}-e^{-\frac{h}{2}H}}{e^{\frac{h}{2}}-e^{-\frac{h}{2}}}
\end{equation}
If $K=e^{\frac{h}{4}H}$ we have the relations

$KX=e^{\frac{h}{2}}XK$

$KY=e^{-\frac{h}{2}}YK$
\begin{equation}
[X,Y]=\frac{K^2-K^{-2}}{e^{\frac{h}{2}}-e^{-\frac{h}{2}}}.
\end{equation}

Notice that there is a Hopf algebra structure on $U_h$ as a module over
 ${\bf C}[[h]]$.

Following \cite{Re-Tu} we consider  $A$ be the quotient of $U_h$ obtained by
 setting

$h=\frac{2\pi\sqrt{-1}}{r}$, $X^r=Y^r=0$, $K^{4r}=1$.

Then $A$ is a finite dimensional algebra over the complex numbers with
 generators
$X,Y,K,K^{-1}$ and the relations stated above. As in the case of $U$ there are
unique
$A$-modules $V^k$ in each dimension $k$ but $V^k$ is irreducible only if
 $k\leq r$.
Also $A$ acquires a Hopf algebra structure from $U_h$ and so tensor products
 and duals
of $A$-modules are still $A$-modules. Moreover the following Clebsch-Gordon
rules
remain valid
\begin{equation}
V^k\otimes V^l=\bigoplus_{p=\mid l-k\mid+1; p+k+l={\small odd}}V^p,
\mbox{ if } k+l\leq r+1
\end{equation}
This $A$ is a {\it quasi-triangular} Hopf algebra (see \cite{Drin}):
there exists an invertible element
$R\in A\otimes A$ satisfying

$R\Delta(u)R^{-1}=\check{\Delta}(u)$, $u\in A$,

$(\Delta\otimes 1)(R)=R_{13}R_{23}$,

$(1\otimes \Delta)(R)=R_{13}R_{12}$, \\
where $\check{\Delta}=P\Delta$, $P$ is the permutation  endomorphism of
$A\otimes A$,
$P(u\otimes v)=v\otimes u$, $R_{12}=R\otimes 1$, $R_{23}=1\otimes R$,
$R_{13}=(P\otimes 1)R_{23}$.

Specifically $R$ may be given by
\begin{equation}
R=\frac{1}{4r}\sum_{0\leq n,a,b}^{n<r;a,b<4r}\frac{e^{\frac{h}{2}}-
e^{-\frac{h}{2}}}{[n!]}t^{ab+(b-a+1)n}\left(X^nK^a+Y^nK^b\right)
\end{equation}
where $t=e^{-\frac{2\pi\sqrt{-1}}{4r}}$.

Remark that $R$ may be viewed as acting on tensor products of two
 $A$-modules $V\otimes W$. We set
$\check{R}:V\otimes W\longrightarrow W\otimes V$ be the {\it flip R-matrix}
$\check{R}=P\circ R$.

(3.3) {\bf Colored framed tangle operators}. Assume the quasi-triangular
Hopf algebra $A$ is fixed.
A coloring of a tangle $T$ is the assignment of an $A$-module to each
of its components.
This way a coloring of $\partial T$ is induced: if $s$ is an arc colored
by $V$ then
assign $V$ to the endpoint of $s$ where is oriented down, and the module
$V^*$ to the
other one. Tensoring from left to the right the modules associated to
the bottom (or upper)
endpoints we get the boundary $A$-modules assigned to $\partial_-T$ and
$\partial_+T$, which we denote $T_-$ and $T_+$ respectively. By convention
the empty product
is ${\bf C}$.

We have two composition laws on tangles $\circ$ and $\otimes$ illustrated in
picture 4.
\begin{figure}
\centering
\hspace{1,5cm}
\psfig{figure=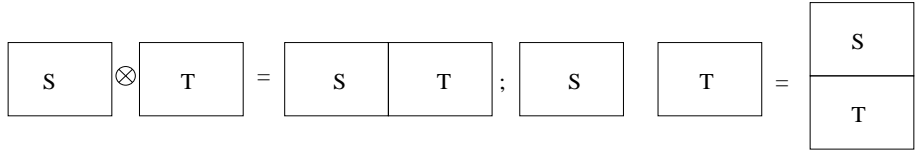,width=12cm}
\caption{Composition laws for tangles}
\end{figure}
\newtheorem{blb}{Theorem}[section]
\begin{blb}
(\cite{Re-Tu,Ki-Me}) There exist uniquely A-linear operators
$J_T:T_-\longrightarrow T_+$ assigned to each colored tangle which satisfy

$J_{S\circ T}=J_S\circ J_T$,

$J_{S\otimes T}=J_S\otimes J_T$,

and for elementary tangles are defined by

$J_I=1$, $J_R=\check{R}$, $J_L=\check{R}^{-1}$, $J_{CR}=E$,
 $J_{CL}=\check{E}$, $J{AR}=N$, $J_{AL}=\check{N}$, \\
where
$R$ is the right hand twist (the orientation points down) tangle, $L$ is the
left
 hand twist,
$CR$ (respectively $AR$) is the creation (annihilation) tangle with the sense
of
 the orientation from left to the right,
$CL$ (respectively $AL$) have opposite orientation than $CR$ and $AR$
respectively,

$E(f\otimes x)=f(x)$, $\check{E}(x\otimes f)=f(K^2x)$,

$N(1)=\sum_{i}e_i\otimes e^i$, $\{e_i\}$ is an arbitrary basis and $\{e^i\}$
its dual,

$\check{N}(1)=\sum_ie^i\otimes K^{-2}e_i$.
\end{blb}
Notice that $J_K$ is just a scalar if $K$ is a colored link.

Now we restrict ourselves to colorings by irreducible $A$-modules so the
set of colors correspond to $\{1,2,...,r\}$. We denote by ${\bf k}$ the
coloring of
$T$ where the $j^{th}$ component is colored with the module of dimension
$k_j$. Thus the theorem yields a family of topological invariants $J_{T,{\bf
k}}$
for colored tangles.

Remark that $J_{T,{\bf k}}$ are independent on the various orientations of
closed
components
of $T$. Also from (\cite{Ki-Me} p.506) if a color in the vector ${\bf k}$ is
$r$ then
the invariant $J_{T,{\bf k}}=0$. So we may assume the colors are from the
subset
$\{1,2,...,r-1\}$.

(3.4) {\bf Closed 3-manifold invariants}. Let $L$ be a framed link in $S^3$.
Recall $L$ determines a 4-manifold $W_L$ obtained by adding 2-handles to the
4-ball $B^4$ along the components of $L$ in $S^3=\partial B^4$. The manifold
$D(L)=\partial W_L$ oriented "outward first" is the result of Dehn surgery on
$L$,
 and
any 3-manifold may be obtained this way. We can pass from one surgery link $L$
 for
$M=D(L)$ to another link $L'$ with $D(L')=M$ by a finite sequence of Kirby
moves
(blow-ups and handle slidings) or equivalently $m$-strands K-moves (see
\cite{Ki-Me}).

Define for a framed link $L$
\begin{equation}
Z_L=\alpha_L\sum_{{\small coloring} {\bf k}}[k]J_{L,{\bf k}},
\end{equation}
where $\alpha_L=b^{n_l}c^{\sigma(L)}$, $[{\bf k}]=\prod [k_i]$,
$[k]=\frac{e^{\frac{kh}{2}}-e^{-\frac{kh}{2}}}{e^{\frac{h}{2}}-
e^{-\frac{h}{2}}}=\frac{\sin\frac{\pi k}{r}}{\sin\frac{\pi}{r}}$,
$b=\sqrt{\frac{2}{r}}\sin\frac{\pi}{r}$,
$c=e^{-\frac{6\pi\sqrt{-1}(r-2)}{8r}}$,
$n_L$ is the number of components of $L$, $\sigma_L$ is the signature of the
 linking matrix
of $L$.

The invariance of $Z_L$ to Kirby moves is proved in \cite{Re-Tu}. Now if $M$ is
obtained
by Dehn surgery on the framed link $L$ we set
$Z_r(M)=Z_L\in {\bf C}$ which is a topological invariant for 3-manifolds.

(3.5) {\bf Cobordisms and TQFT}. For simplicity we restrict ourselves to
cobordisms
$M$ with connected boundaries $\partial_-M$ (incoming) and $\partial_+M$
(outgoing).

We define first the spaces associated to closed oriented surfaces. Consider
$(\Gamma, {\bf i})$ be one of the
colored 3-valent graph of genus $g$  from picture 5, viewed as a
(0,2g)-tangle in ${\bf R}^3$.
\begin{figure}
\centering
\hspace{1cm}
\psfig{figure=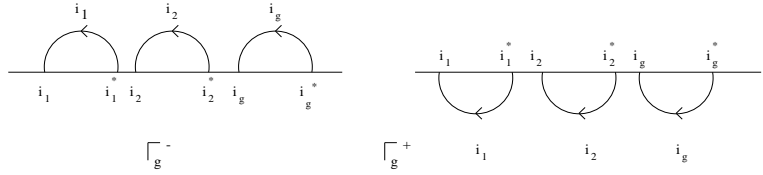,width=10cm}
\caption{The standard spine for $\Sigma_g$}
\end{figure}

We associate to this colored graph the space

$Z(\Gamma_g,{\bf i})=\bigotimes_{l=1}^g(V^{i_l}\otimes V^{i_l*})$. \\
Then to the closed oriented surface of genus $g$ we assign the space

$Z(\Sigma_g)=\bigoplus_{{\small coloring} {\bf k}}Z(\Gamma_g,{\bf k})$.

We may extend now the Dehn surgery construction to cobordisms using 3-valent
graphs. We call $\Gamma$ a special framed graph if it satisfies the conditions:

(i) $\Gamma\cap I^2\times[0,\frac{1}{10}]$ is the base of a standard 3-valent
graph $\Gamma^+_{g_+}$. This means that
$\Gamma\cap I^2\times[0,\frac{1}{10}]$ from which the components not touching
$I^2\times\{0\}$ are removed is isomorphic to  $\Gamma^+_{g_+}$.

(ii) $\Gamma\cap I^2\times[\frac{9}{10},1]$ is the base of the standard
3-valent
graph
$\Gamma^-_{g_-}$.

(iii) $\Gamma\subset I^3$ and the union of the components  which do not touch
the boundaries form a link $L$.

Now each special framed link $\Gamma$ gives rise to a decorated cobordism
$(M,\Sigma_-,\Sigma_+)$ with parametrized surfaces $\Sigma_-, \Sigma_+$ of
genera
$g_-,g_+$ respectively. The construction goes as follows: we have a regular
 neighborhood
$N(\Gamma^+_{g_+})\subset S^3$ and an homeomorphism
$f^+:H_g\longrightarrow N(\Gamma_{g_+}^+)$
from the handlebody of genus $g$ (and a similar situation  for
$\Gamma^-_{g_-}$).
Cut out open handlebodies $N(\Gamma_{g_+}^+)$ and
$N(\Gamma_{g_-}^-)$ from $S^3$ to get a compact oriented 3-dimensional
cobordism $E$ between the respective surfaces. Now the maps $f^+,f^-$ induce
parametrizations of the boundary of $E$.
 Then surgery on $E$ on the remaining link
$L$ produces a compact cobordisms $M=D(\Gamma)$ with parametrized
boundaries $\partial_-M\cong \Sigma_{g_-}$ and
$\partial_+M\cong \Sigma_{g_+}$. Again each cobordism $M$ whose boundary is
partitioned into two disjoint parts may be obtained this way, and there are
generalized Kirby moves for such surgery
presentations (see \cite{Tu2} p.168, \cite{Fun5}).

We set now
\begin{equation}
Z^{{\bf i}}_{{\bf j}}(\Gamma)=b^{n_L-g_+}c^{\sigma(L)}[j]
\sum_{{\small coloring} {\bf k}}[{\bf k}]J_{\Gamma,({\bf i,j,k}}
\end{equation}
where ${\bf i},{\bf j}$ are colorings of $\Gamma_{g_-}^-,\Gamma_{g_+}^+$
respectively,
${\bf k}$ is coloring of $L$, and $(\Gamma,{\bf i,j,k})$ is viewed as a
colored
tangle giving rise to a map
$J_{\Gamma,{\bf i,j,k}}:Z(\Gamma_{g_-}^-,{\bf i})\longrightarrow
Z(\Gamma_{g_+}^+,{\bf j})$.
We set $Z(\Gamma):Z(\Sigma_{g_-})\longrightarrow Z(\Sigma_{g_+})$ for the
 linear map
whose blocks are the matrices $Z^{{\bf i}}_{{\bf j}}(\Gamma)$.

As $Z(\Gamma)$ is invaried by Kirby moves it follows that the formula
$Z_r(M)=Z(\Gamma)$ defines a topological invariant for 3-dimensional oriented
cobordisms with parametrized boundaries (i.e. decorated cobordisms).

Moreover we have

$Z_r(M\circ N)=c Z_r(M)\circ Z_r(N)$, \\
where $c$ lies in the group of roots of unity generated by
 $e^{\frac{2\pi\sqrt{-1}}{r+2}}$.
Therefore this data is a TQFT with anomaly which we call the
$sl_2({\bf C})$-TQFT at level $r$.

\section{The proof of the theorem}

(4.1) {\bf The Arf invariant}. We have a recurrent method to compute the
Arf invariant (see \cite{Robert}) of a proper link due to Murakami \cite{Mura}.
This is related to Jones polynomial at 4-th roots of unity. Specifically let
$I$ denote the link invariant defined by

$I(unknot)=1$

$I(unknot \bigcup K)=\sqrt{2} I(K)$, for any link $K$, \\
and the skein relation

$I(L_+)+I(L_-)=\sqrt{2} I(L_0)$, \\
where $L_+,L_-$ are the left and right hand twists and $L_0$ is the 2-parallel
string diagram, and the rest of the
diagrams are the same (see picture 6).
\begin{figure}
\centering
\hspace{1cm}
\psfig{figure=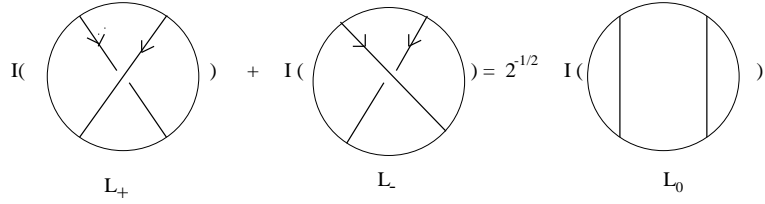,width=10cm}
\caption{The skein relation}
\end{figure}

Then for a link $L$ we have

$I(L)=\left\{\begin{array}{lll}
       (-1)^{\varepsilon}\sqrt{2^{n_L-1}} & \mbox{ if } & Arf(L)=\varepsilon \\
       0 & \mbox{ if } & L \mbox { is not proper}
       \end{array}\right.$\\
Remember that the link $L$ is proper if, for each sub-link $K$ the linking
number
$lk(K,L-K)$ is even.

(4.2) {\bf Jones polynomial at 4th roots of unity}.
Consider now $J_L=J_{L,{\bf 2}}$, where ${\bf 2}$ is the coloring of all
components of $L$ by
the module $V^2$. Then $J_L$ is a variant of Jones polynomial according to
\cite{Ki-Me}. In fact we have :

$J_L=t^{3L\cdot L}\sqrt{2}I(L)$, \\
where $t=e^{\frac{2\pi \sqrt{-1}}{16}}$.

4.3) {\bf The cabling formula}. If $L$ is a framed link ${\bf k}$ a coloring
then
the following formula permits to compute $J_{L,{\bf k}}$ in terms only of
Jones polynomial of cablings of $L$. Specifically we have
\begin{equation}
J_{L,{\bf k}}=\sum_{{\bf j}=
{\bf 0}}^{{\bf n}/2}(-1)^{\bf j}C_{{\bf n}-{\bf j}}^{\bf j}J_{L^{{\bf n}-2{\bf
j}}}
\end{equation}
where ${\bf n}={\bf k}-{\bf 1}$, and we set $f({\bf n})=\prod_if(n_i)$,
${\bf m}<{\bf n}$ if
$m_i<n_i$ for all $i$'s etc. Also $L^{\bf c}$ is the ${\bf c}$-cabling of $L$
which consists in replacing the $i^{th}$ component of $L$ by $c_i$ parallel
copies.

In case when $r=4$ and the link $L$ has two components $K$ and $H$, both
unknotted in ${\bf R}^3$, then the possible values for $J_{L,{\bf k}}$ are

$J_{L,(1,1)}=1$,

$J_{L,(1,2)}=J_H=\sqrt{2}$,

$J_{L,(2,1)}=J_K=\sqrt{2}$,\\
and using the 1-colored components removing lemma (see \cite{Ki-Me},p.511),

$J_{L,(1,3)}=J_{H,3}=J_{H^2}-1=1$,

$J_{L,(3,1)}=J_{K,3}=J_{K^2}-1=1$,

$J_{L,(2,2)}=J_L$

$J_{L,(3,2)}=J_{K^2H}-J_H=J_{K^2H}-\sqrt{2}$,

$J_{L,(2,3)}=J_{KH^2}-J_K=J_{KH^2}-\sqrt{2}$,

$J_{L,(3,3)}=J_{K^2H^2}-J_{K^2}-J_{H^2}+1=J_{K^2H^2}-3$.

(4.3) {\bf A surgical description of the cobordism $X$}.
We come back now to the cobordism $X$ from (2.4).
We have a simple surgical description for $X$ since both tori $T_0$ and $T_1$
are unknotted in $S^3$. We can choose for example the special graph $\Gamma$
from
picture 7.
\begin{figure}
\centering
\hspace{1cm}
\psfig{figure=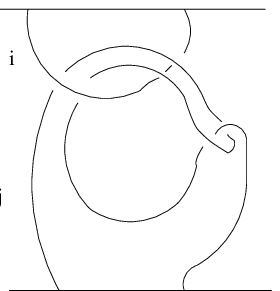,width=3cm}
\caption{The special graph $\Gamma$}
\end{figure}
Taking into account that the intermediary link of $\Gamma$ is trivial this time
we see that

$Z^i_j(\Gamma)=c^{-3}J_{\Gamma,(i,j)}$, for $i,j\in\{1,2,3\}$. \\
 When properly interpretated
$J_ {\Gamma,(i,j)}$ is $J_{L,(i,j)}$ where $L$ is the Whitehead link (see the
 picture 8).
\begin{figure}
\centering
\hspace{1cm}
\psfig{figure=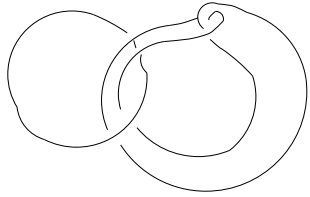,width=4cm}
\caption{The Whitehead link}
\end{figure}
Since both components of $L$ are unknotted we may apply the previous
formulas.
We compute first

$I(L)=-\sqrt{2}$, $I(K^2H)=I(KH^2)=-2$, $I(K^2H^2)=4$.

Therefore up to a root of unity the morphism $Z(X)$ is given by the matrix

$\centering
\left(\begin{array}{ccc}
       1 & \sqrt{2} & 1 \\
       \sqrt{2} & \sqrt{2}(1-\sqrt{-1}) & 2\sqrt{-2}-\sqrt{2} \\
       1 & 2-2\sqrt{-1}-\sqrt{2} & -3-4\sqrt{-1}
       \end{array}\right)$ \\
whose determinant is $2-13\sqrt{2}+(16-4\sqrt{2})\sqrt{-1}$. Therefore
the inductive limit $\lim_{\rightarrow}({\bf C},Z(X))$ of
iterates of the map $Z(X)$ is isomorphic to ${\bf C}^3$ and the claim of the
theorem
is proved. Notice that the change of the parametrization on intermediary
 boundaries
amounts to multiply the matrix $Z(X)$ by an invertible one which does not
affect
the limit.

\bibliography{ref}

\bibliographystyle{plain}

\end{document}